\documentclass{svproc}
\usepackage{hyperref}
\usepackage{url}
\usepackage{amsmath}
\usepackage{graphicx} 
\usepackage{float}

\begin{document}
\mainmatter              
\title{Solving the option Forecast problem by a numerical method for the Black-Scholes Equation with Machine Learning Classification Model}
\titlerunning{Classification Model for Option Forecasting}  
%
\author{Benjamin Jiang\inst{1} \and Matthieu Durieux\inst{2} \and 
Kirill V. Golubnichiy\inst{3}}
\authorrunning{B. Jiang et al.} 
%
\tocauthor{First Author, Second Author, Third Author}
\institute{University of Washington, Seattle WA 98195, USA,\\
\email{jianghao@uw.edu}
\and
University of Washington, Seattle WA 98195, USA,\\
\email{matt45@uw.edu}
\and
Texas Tech University, Lubbock TX 79409, USA,\\
\email{kgolubni@ttu.edu}}

\maketitle              

\begin{abstract}
In this study, we propose novel classification models that integrate the Quasi-Reversibility Method (QRM) with advanced machine learning techniques to enhance the prediction of option prices. The QRM, well known for solving the Black-Scholes equation under challenging conditions, forecasts option prices one day in advance \cite{KlibGol}. By leveraging the QRM-derived minimizer with machine learning classification models, we aim to categorize options as either increasing or decreasing in value. This integration enables more refined categorization of financial instruments, blending numerical analysis and predictive modeling.

Using these classifications, we implement trading strategies to enhance the predictive power of QRM extrapolations. These strategies leverage both historical data trends and QRM forecasts, potentially improving trading efficiency, especially in volatile markets. To validate our approach, we collected 23,548 real-world options data points across various time periods and market conditions, creating a comprehensive dataset to rigorously test our models. This data, combined with the QRM minimizer, trains machine learning models, including decision trees, random forests, gradient boosting classifiers, k-nearest neighbors, and neural networks.

The models’ performance is assessed using metrics such as accuracy, precision, and recall to evaluate the framework’s consistency and reliability. Our goal is to identify which models, paired with QRM-derived features, perform most effectively for option price prediction. This study ultimately aims to provide a robust predictive tool that is practical for short-term trading and grounded in established mathematical finance methods.
\end{abstract}

\section{Introduction and Overview}
The stock market, also known as the equity market, consists of buyers and sellers trading ownership of businesses. To profit in such a market, one must intuitively guess market trends, aiming to buy low and sell high. However, due to the market's volatility and the ill-posed nature of predicting exact stock prices caused by noisy data and model sensitivity, accurate predictions remain challenging. Therefore, a more generalized approach, such as estimating market gains or losses, is often more feasible.

Options trading involves contracts that grant the right to buy or sell stock in the future. Predicting whether an option will trade higher can lead to market profits. We developed the Quasi-Reversibility Method (QRM) to capture market movements by solving the ill-posed problem of the Black–Scholes model to predict option prices one day in advance \cite{KlibGol}. Our research utilizes the QRM-derived minimizer, which predicts future option data, combined with other stock market features to forecast market movements. By integrating these option price predictions, we aim to improve the accuracy of predicting whether the overall market will gap up or down, acknowledging the short-term inefficiencies that may exist in the market.

We propose using ensemble methods, such as gradient boosting trees and random forests, as well as other advanced classification techniques including k-nearest neighbors, and neural networks, to classify future market data. By combining these classification models with QRM, we aim to capture critical market features that influence overall price movements, assuming limited market volatility in the short term. Our approach involves designing weak learners to classify market movements into different categories and then improving the results by minimizing the loss function.

This paper provides an introduction to the Quasi-Reversibility Model (QRM) employed in our research, along with a brief summary of prior work in this area. Additionally, we present various classification models, including gradient boosting, random forests, k-nearest neighbors, and neural networks, and explain how these models are applied using the results obtained from QRM. The findings are demonstrated through an analysis of key features and confusion matrices. Following this, we provide a comprehensive summary of critical evaluation metrics such as precision, recall, and accuracy. Finally, we discuss future developments and potential directions for this research. If the proposed classification models yield favorable results, they could enhance the profitability of options trading by improving the prediction of market movements.
\section{Quasi-reversibility Model and Black–Scholes Model}
\subsection{Black-Scholes Equation}
The Black-Scholes equation is a backward parabolic partial differential equation used in mathematical finance. The proposed model uses initial boundary conditions for the stock of interest to simulate a more accurate options trading strategy. Since it is an ill-posed problem, we introduce the variable $\tau$ to solve the equation backward, where at any given time $t$, we have:
\begin{align}
\tau = T - t \tag{1}
\end{align}
The Black-Scholes equation is given by:
\begin{align}
\frac{\partial{V}}{\partial{t}} + \frac{1}{2}\sigma^2 S^2 \frac{\partial^2{V}}{\partial{S^2}} + rS\frac{\partial{V}}{\partial{S}} - rV = 0 \tag{2}
\end{align}
Assuming that the function $u(s, \tau)$ satisfies the Black-Scholes equation with a volatility coefficient $\sigma$, we can rewrite the equation as:
\begin{align}
\frac{\partial{u(x, \tau)}}{\partial{\tau}} + \frac{1}{2}\sigma^2 x^2 \frac{\partial^2{u(x, \tau)}}{\partial{x^2}} = 0 \tag{3}
\end{align}
with the initial condition:
\begin{align}
u(x,0) = f(x) \tag{4}
\end{align}
where $f(x)$ is the payoff function, $x$ is the stock price, and at maturity time $\tau = 0$. We assume the risk-free interest rate to be zero, and thus, $u(x, \tau)$ represents the price of an option at time $t$ \cite{function}.
\subsection{Quasi-Reversibility Method}
We introduce this method to solve Black-Scholes equations with certain approximations. We aim to solve the following equation:
\begin{align}
Pu &= u_t + \frac{1}{2}\sigma^2(t)x^2u_{xx}=0 \tag{5}
\end{align}
with the Dirichlet boundary conditions:
\begin{align}
u(x_b, t) = o_b(t), \quad u(x_a, t) = o_a(t) \tag{6}
\end{align}
and initial conditions:
\begin{align}
u(x, 0) = f(x) \tag{7}
\end{align}
where:
\begin{align}
f_x &= \frac{x_a}{x_b} - 1 \tag{8}
\end{align}
\begin{align}
f_o &= \frac{o_a}{o_b} - 1 \tag{9}
\end{align}
with $t \in [0,2\tau]$, and $x \in (x_b(0),x_a(0))$. $P$ represents the partial differential operator for the Black-Scholes equation, where $x$ is the stock price, $\sigma$ is the volatility of the option, $o_b$ and $o_a$ are the bid and ask prices of the option, and $x_b$ and $x_a$ are the bid and ask prices of the stock. We assume that the bid price is always less than the ask price.

The method involves five steps to solve the Black-Scholes equations:
\begin{itemize}
\item First, solve for the dimensionless variable, namely $x$.
\item Secondly, perform interpolation and extrapolation to forecast the option price.
\item Thirdly, apply the initial boundary conditions to the PDE: Partial Differential Equation.
\item Fourthly, apply regularization to solve the method numerically.
\item Lastly, solve the minimizing problem using finite difference and regularization \cite{function}.
\end{itemize}
\subsection{Analysis}
Assuming $Y_{2\tau} = (x_b(0), x_a(0)) \times (0, 2\tau)$, the functions:
\begin{align}
o_b(t), o_a(t) \in H^2[0, 2\tau], \quad \sigma(t) \in C^1[0, 2\tau] \tag{10}
\end{align}
allow us to find the solution $u \in H^2(Y_{2\tau})$ that satisfies the following conditions:
\begin{align}
Au = 0 \text{ in } Y_{2\tau} \tag{11}
\end{align}
\begin{align}
u(0, t) = o_b(t), \quad u(1, t) = o_a(t) \tag{12}
\end{align}
\begin{align}
u(x,0) = f(x), \quad x \in (0,1) \tag{13}
\end{align}
We will perform a general convergence analysis based on a new Carleman estimate, assuming the time interval is within a short period. Let's consider a general case with the given properties. Given a non-negative constant $L$ and let $Z_L$ be $(x,t) \in (0,1) \times (0, L)$. Suppose there are two non-negative constants $\gamma_1, \gamma_2$ with $\gamma_2 > \gamma_1$. Suppose the functions $\Gamma(x,t) \in C^1(\bar{Z}_{L})$ satisfy:
\begin{align}
||b||_{C^1(\bar{Z}_{L})} \leq \gamma_2, \quad \Gamma(x,t) \geq \gamma_1 \text{ in } Z_L \tag{14}
\end{align}
Then, the general version of the original problem is formulated as:
\begin{align}
Mk = k_t + \Gamma(x,t)k{xx} = 0 \text{ in } Z_L \tag{15}
\end{align}
with the boundary condition:
\begin{align}
k(0, t) = o_b(t), \quad k(1, t) = o_a(t) \tag{16}
\end{align}
and the initial condition:
\begin{align}
k(x,0) = k(x) = o_b(0)(1-x) + o_a(0)x, \quad x \in (0, 1). \tag{17}
\end{align}
We are trying to find a convergent solution to the general solutions \cite{SH}.
Our previous initial results using neural networks (NN) were presented in the paper \cite{function}, where we tested three different methods: the quasi-reversibility method for solving ill-posed problems, and two models using multilayer fully connected NNs—classification and regression. 
\begin{table}
\caption{Results of the previous model}
\label{table:prev_results}
\begin{center}
\begin{tabular}{c| c| c| c} 
 \hline
 Method & Accuracy & Precision & Recall\\ 
 \hline
 QRM & 49.77\% & 55.77\% & 52.43\% \\
 \hline
 Binary Classification & 56.36\% & 59.56\% & 70.22\% \\
  \hline
 Neural Network & 55.42\% & 60.32\% & 61.29\% \\ 
 \hline
\end{tabular}
\end{center}
\end{table}

\section{Classification Algorithms}
\subsection{Introduction}
Classification algorithms are a subset of supervised machine learning techniques that aim to categorize input data into predefined classes or labels. These algorithms learn from labeled training data, where the input features are associated with known output labels, to build a model that can predict the class of new, unseen data. Common classification algorithms include decision trees, random forests, k-nearest neighbors (KNN), and neural networks. Each algorithm has its strengths and is chosen based on the nature of the dataset and the specific requirements of the problem. For instance, decision trees are easy to interpret and visualize, while neural networks can model complex relationships in data. The performance of classification algorithms is typically evaluated using metrics such as accuracy, precision, recall, and the F1 score. These metrics provide insights into how well the model distinguishes between different classes. Properly selecting and tuning these algorithms is crucial for developing effective predictive models in various domains, from medical diagnosis to stock market forecasting \cite{scikits}.
\subsection{Decision Tree}
A decision tree is a classification algorithm that utilizes a tree-like model. This model includes a root node with no incoming edges, decision nodes with exactly one incoming edge and at least one outgoing edge, and leaf nodes with only incoming edges and no outgoing edges. The incoming edge represents the result of a previous decision, while the outgoing edges represent possible outcomes of current decisions. It is also possible to have branches consisting of multiple nodes. Note that the decision tree classifies data into final categories based on these decision paths \cite{Tree}. 

\begin{figure}
\centering
\includegraphics[width=0.68\textwidth]{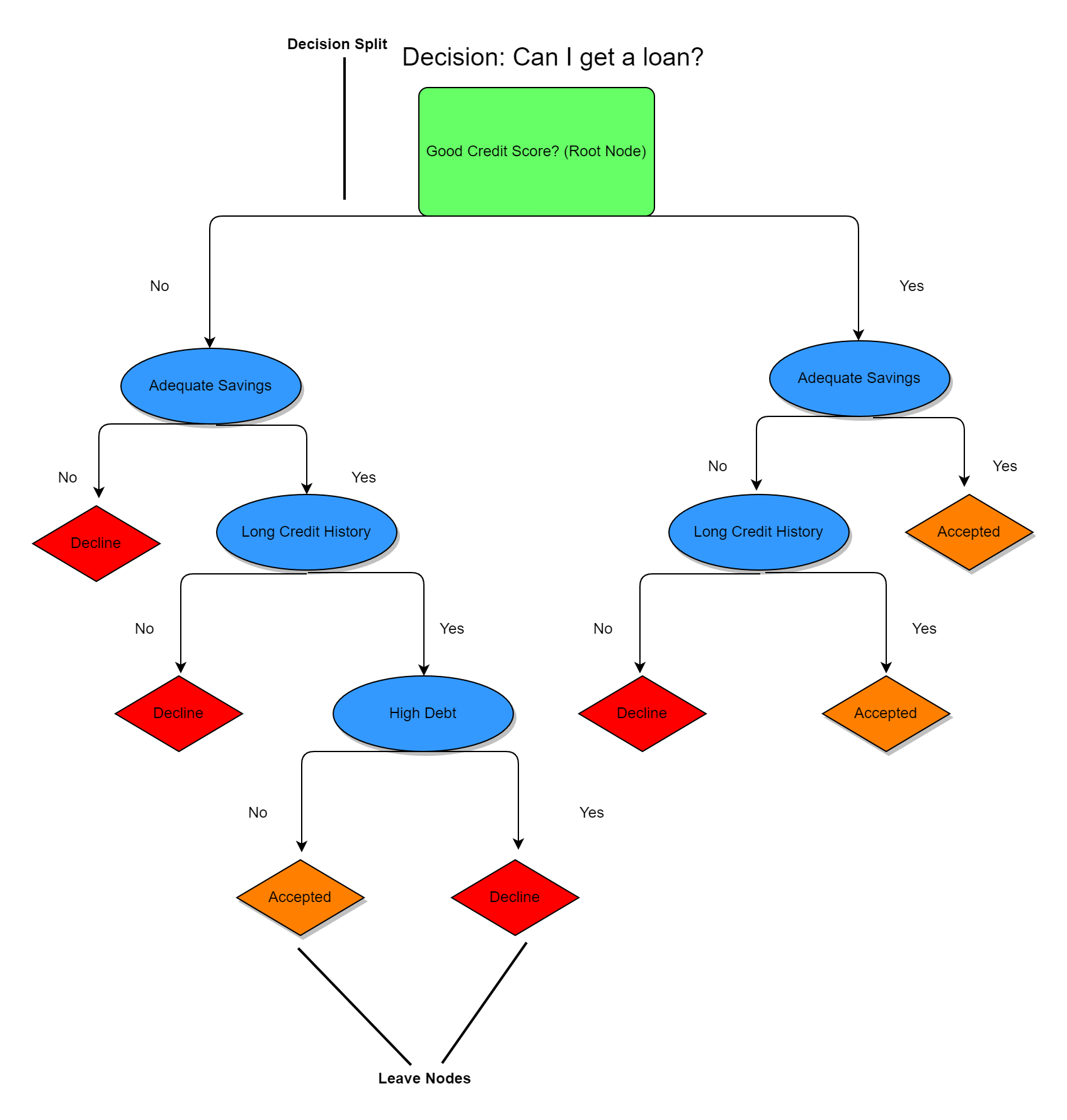}
\caption{Binary Decision Tree example}
\label{fig}
\end{figure}

The concept behind a decision tree is straightforward. By using specific features of the data, the tree splits the data at each node based on the values of these features, with each outgoing edge from a node representing a possible outcome from the decision made. While a decision tree does not have to be binary—meaning nodes can have more than two outgoing edges—we will use a binary tree in our case because we want to capture the binary movement of options as either increasing or decreasing. We will feed the tree with data derived from the QRM minimizer and market features, splitting the data into training, validation, and test sets. Subsequently, we will create trees with different combinations of features to generate outcomes for our trading strategy. Below is an illustration showing a basic binary decision tree, which we refer to as a classification tree (see Figure \ref{fig}).

\subsection{Data Split}
Data splitting is the process used to determine the decisions at each node of a decision tree. Generally, we select a feature that divides the data roughly in half to serve as the root node, which helps to keep the tree balanced. At the same time, we design the leaf nodes to have a bias towards one of the classifying classes, enabling us to make a final decision.
\subsubsection{Entropy}
Entropy is a measure of the disorder or unpredictability of a feature relative to our data. It determines how effectively the data can be classified into a specific class. The formula for entropy is given below:
\begin{align}
H(S) = \sum_{i=1}^{c}-p_i \log_{2}{p_i} \tag{18}
\end{align}
where $H(S)$ represents the entropy or disorder of the feature, $p_i$ is the probability of encountering feature $i$, and $c$ is the total number of features. Entropy can range from 0 to 1, with 0 indicating the lowest entropy, signaling no disorder. This is undesirable for training as it implies that all data points have been classified into a single class. Conversely, an entropy of 1 indicates a balanced data split, which is ideal for training. Moreover, entropy is used to calculate the information gain at each node.
\subsubsection{Information Gain}
Information Gain measures the change in entropy resulting from splitting the data on a feature. The formula is as follows:
\begin{align}
\mbox{Information Gain}(S, X) = H(S) - \sum_{i=1} \frac{x_i}{S} H(x_i) \tag{19}
\end{align}
where $H(S)$ represents the measure of disorder before the split, $x_i$ denotes the number of observations in the $i^{th}$ split, and $H(x_i)$ is the disorder of the target feature in split $x_i$. By subtracting the combined entropy of all child nodes from the entropy of the parent node, we can determine how much "information" we have gained or lost. A higher Information Gain indicates a larger reduction in entropy or disorder. Information Gain is used to assess whether a split was beneficial by evaluating how much "noise" was eliminated due to the split.
\subsubsection{Gini Impurity}
Gini Impurity measures the probability that an instance of a random variable is incorrectly classified if that instance were randomly classified according to the distribution of class labels from the dataset. The formula is presented below:
\begin{align}
Gini = 1 - \sum_{i=1}^{n}(p_i)^2 \tag{20}
\end{align}
where $n$ is the number of classes, and $p_i$ is the probability of an item being labeled as class $i$. To calculate Gini Impurity, one sums the squared probabilities for each class, $p_i^2$, and subtracts this value from 1. Gini Impurity varies between 0 and 1, where 0 represents complete purity (all elements belong to a single class), 0.5 represents an equal distribution of elements across classes, and 1 indicates maximum impurity, with elements uniformly distributed across multiple classes.
\subsection{Hyperparameter Selection}
Hyperparameters define a model much like how the rules of a sport define that sport. One can differentiate between models through their hyperparameters, as these address questions related to model design. Such questions might include determining the maximum depth of a tree, the maximum number of nodes, or the learning rate for gradient descent. Hyperparameters are not features of the model but rather constitute its structure. They determine how the algorithm behaves during execution, making tuning these hyperparameters essential to minimize noise and enhance accuracy. The selection process can vary; one might manually adjust hyperparameters hoping to achieve better results, or use a method such as Grid Search. Grid Search involves building a model for each possible combination of hyperparameters, evaluating each model, and selecting the one with the most promising results. Another method, Random Search, generates random hyperparameter values from a data distribution for each parameter. This process is called Random Search due to the randomness introduced by not using a predefined set of values, but rather randomly sampling the values from a distribution \cite{Skit}.
\subsection{Ensemble Methods}
The ensemble method utilizes groups of base learners or models to form a final decision by averaging the outputs of different learners or models. This approach is akin to voting, where the decision is based on the majority opinion. Since individual models might exhibit high bias or variance, combining models can reduce these issues by either taking random samples from the same model or using different models to train the data. This principle underlies the bagging and boosting methods, where bagging trains weaker learners in parallel and boosting trains them sequentially. We employ these methods to mitigate the potential for high bias.
\subsection{Gradient Descent}
\begin{figure}
\centering
\includegraphics[width=1\textwidth]{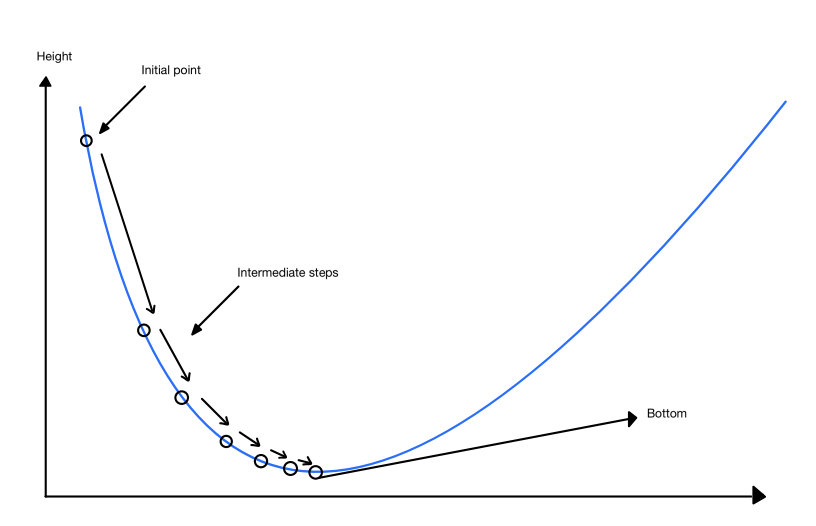}
\caption{Example of Gradient descent for convergence}
\label{fig:descent}
\end{figure}

Gradient descent is a first-order iterative optimization method that seeks to find the local minima of an algorithm. It is a boosting method that aims to reduce inaccuracies in our model by transforming a weak model into a stronger one by addressing its shortcomings. The concept involves observing the gradient at the current position and then taking a step in the opposite direction, which results in the steepest descent towards the minimum. This method is analogous to descending a hill along the steepest path. Additionally, there is a parameter known as the learning rate, which determines the size of each step taken during the descent. A larger learning rate might result in overshooting the minimum and ascending the opposite slope, while a smaller step ensures greater accuracy but may be inefficient and susceptible to saddle points, which can mislead the algorithm into thinking a local minimum has been reached (\ref{fig:descent}. In this method, we use a loss function to gauge our progress down the hill; essentially, it represents the difference between the actual output of our model and our prediction, which we refer to as the error. We strive to minimize the loss function until it approaches zero, ideally reaching the global minimum. However, there are challenges associated with using this algorithm. It is a greedy method that focuses solely on the current location, potentially leading to a local minimum that may not represent the optimal solution \cite{Descent}.

\subsection{Gradient Classification}
Here, we aim to apply gradient descent to our decision tree to enhance accuracy.
First, we introduce the mean squared error, which evaluates the quality of our splits. Since we are addressing an optimization problem—specifically, aiming for a classification to determine the labels of different classes—we will utilize the log-loss function and strive to minimize it for improved accuracy.

The idea is to use some input data along with a differentiable loss function. For classification, we employ the log-loss function, which is defined as follows:
\begin{align}
\mbox{total loss} = -\frac{1}{N}\sum_{i=1}^{N}y_i\cdot \log(p(y_i)) + (1-y_i) \cdot \log(1-p(y_i)) \tag{21}
\end{align}
The logarithmic loss is a metric that evaluates the prediction probability corresponding to the true label.
With this loss function, we can proceed with the iterations:

\begin{enumerate}
\item Initialize our model with a constant value:
\begin{align}
F_0(x) = \arg\min_\gamma \sum_{i=1}^{n}L(y_i, \gamma) \tag{22}
\end{align}
\item Compute the pseudo-residuals for each data point:
\begin{align}
r_{im} = - \left[\frac{\partial L(y_i,F(x_i))}{\partial F(x_i)}\right]{F(x)=F{m-1}(x)}, \forall i = 1,2,...,n \tag{23}
\end{align}
\item Fit a regression tree $h_m(x)$ to the pseudo-residuals.
\item Compute the multiplier $\gamma_m$ by solving the one-dimensional optimization problem:
\begin{align}
\gamma_m = \arg\min_\gamma \sum_{i=1}^n L(y_i, F_{m-1}(x_i) + \gamma h_m(x_i)) \tag{24}
\end{align}
\item Update the model:
\begin{align}
F_m(x) = F_{m-1}(x) + \gamma_m h_m(x) \tag{25}
\end{align}
\end{enumerate}
The final output, $F_m(x)$, will serve as our new decision tree model \cite{Skit}. By differentiating the log-loss function at every iteration, we can significantly improve the accuracy of our decision tree. Additionally, employing cross-validation helps mitigate overfitting in deeper trees that show higher training accuracy.
\subsection{Random Forest}
We employ a Bagging method, which selects a subset of the data to create multiple diverse training sets. These sets are then trained independently, and the majority of the predictions are averaged to obtain a more accurate estimation while reducing overfitting \cite{Forest}. In this project, we intend to use the random forest algorithm, which applies the bagging method to develop numerous random, uncorrelated decision trees—a "forest." This process, known as "the random subspace method," aims to generate trees with low correlation among them. Unlike decision trees that use all provided features, random forests consider only a fraction of the features. The data split is performed based on Gini, entropy, or log-loss criteria, determined through hyperparameter tuning. Additionally, to identify important features, we calculate the importance of a feature by assessing the decrease in node impurity weighted by the probability of reaching that node:

\begin{align}
\text{imp}(i) = w_iN_i - w_{left(i)}N_{left(i)} - w_{right(i)}N_{right(i)} \tag{26}
\end{align}
where:
\begin{itemize}
\item $w_i$ = weighted samples that reach node $i$,
\item $N_i$ = impurity of node $i$,
\item $left$ = left split on node $i$,
\item $right$ = right split on node $i$ \cite{Skit}.
\end{itemize}
Although random forest is known for its accuracy, it can be an extremely time-consuming process because it trains many trees simultaneously. Moreover, interpreting the results from such a large collection of decision trees can be challenging. This is why we utilize feature importance metrics to aid in interpreting the results.
Finally, we choose to apply the random forest algorithm to our data because we believe that the market exhibits certain patterns. By correctly identifying and acting on these patterns, we can enhance the predictions made from the QRM. Furthermore, since our primary goal is to classify future options as either increasing or decreasing, having majority votes from different trees can reduce overfitting and bias within the short period of option price predictions.
\subsection{Neural Networks}
Neural networks are a fundamental component of artificial intelligence, modeled after the human brain's architecture. They consist of layers of interconnected nodes, or neurons, that process input data and adjust their connections to minimize prediction errors through a process known as backpropagation. This iterative learning process allows neural networks to improve their accuracy over time. Neural networks are particularly effective for complex tasks involving large amounts of data, such as image and speech recognition. Advances in deep learning, which involves neural networks with multiple hidden layers, have significantly improved their ability to model intricate patterns and relationships in data \cite{neural}.
\subsection{K-Nearest Neighbors}
The k-nearest neighbors (KNN) algorithm is a simple, yet powerful, machine learning technique used for classification and regression. KNN classifies a new data point based on the majority class of its k-nearest neighbors, determined by a distance metric such as Euclidean distance. As a non-parametric method, KNN does not make assumptions about the underlying data distribution, making it flexible and easy to implement. However, it can be computationally intensive with large datasets due to the necessity of calculating distances for all training samples. KNN is particularly effective for small datasets with clear class boundaries and is frequently employed in pattern recognition and recommendation systems \cite{knn}.
\section{Result and Data}\label{sec:algorithms}
\subsection{Introduction}
The results were generated using a robust combination of machine learning methods, integrating the Quasi-Reversibility Method (QRM) as a significant feature in our analysis. We ensured the integrity and reproducibility of our data by setting the random state to 42. This approach provides consistency across different runs of the experiment. Furthermore, we leveraged real market data sourced from Bloomberg, ensuring that our findings are grounded in actual market conditions. By incorporating diverse methodologies and high-quality data, our results offer a comprehensive and reliable analysis of options pricing and market prediction.
\subsection{Classification Matrix}
We define the option price of a certain day to be the average of the bid and the asking price.
All the results below were calculated based on the following matrix:
\begin{itemize}
    \item Predicted label of Up represents that we predicted the option price would increase the next day.
    \item Predicted label of Down represents that we predicted the option price would decrease the next day.
    \item True label of Up represents that the real option price increases the next day.
    \item True label of Down represents that the real option price decreases the next day.
\end{itemize}
\subsection{Decision Tree}
\begin{figure}[ht!]
    \centering
    \includegraphics[width=1\textwidth]{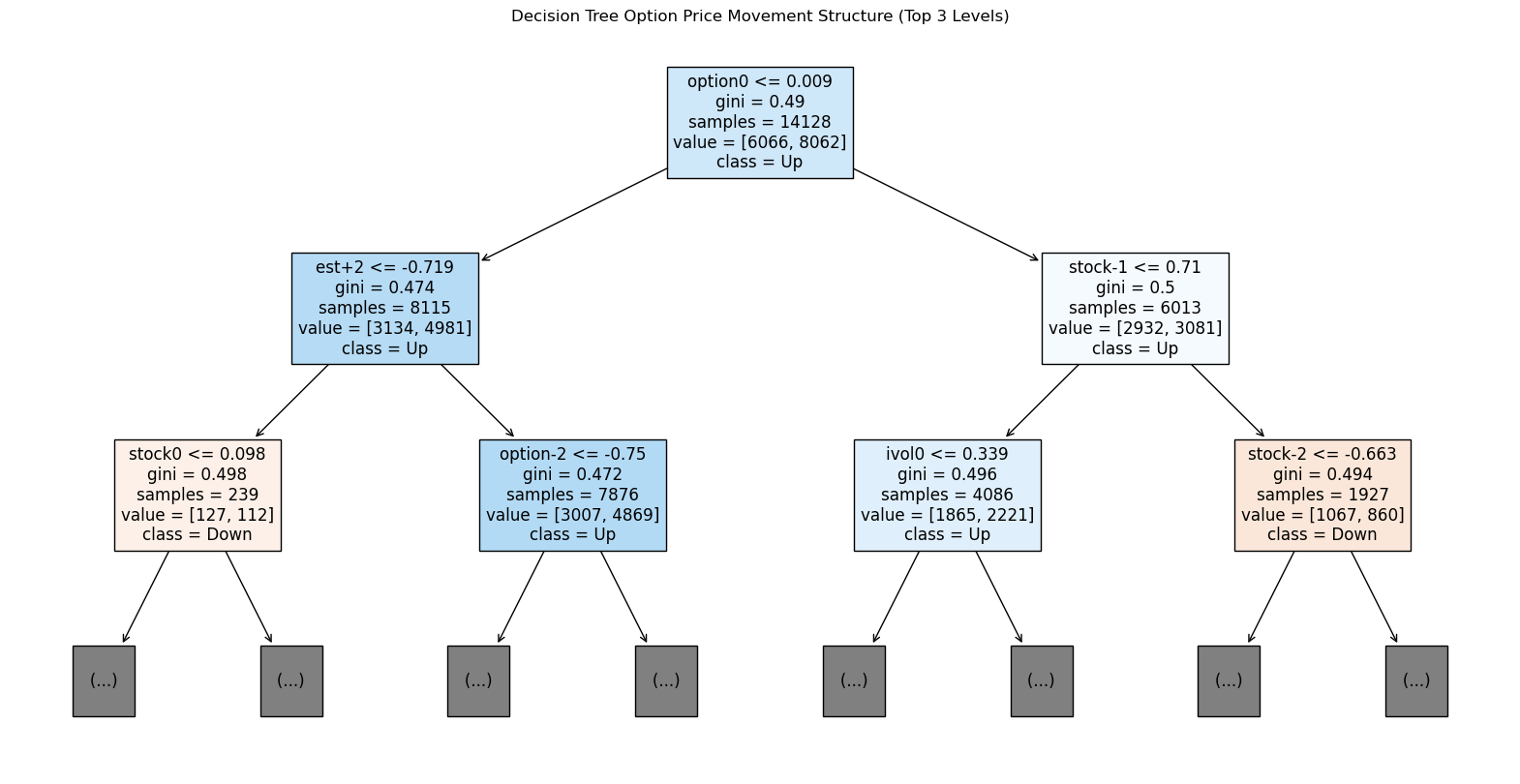}
    \caption{Decision Tree for option movement}
    \label{fig:tree2}
\end{figure}

Since the tree structure is too complex for demonstration purposes, we will only show a partial tree representing the first few layers \ref{fig:tree2}.
\begin{figure}[ht!]
    \centering
    \includegraphics[width=0.7\textwidth]{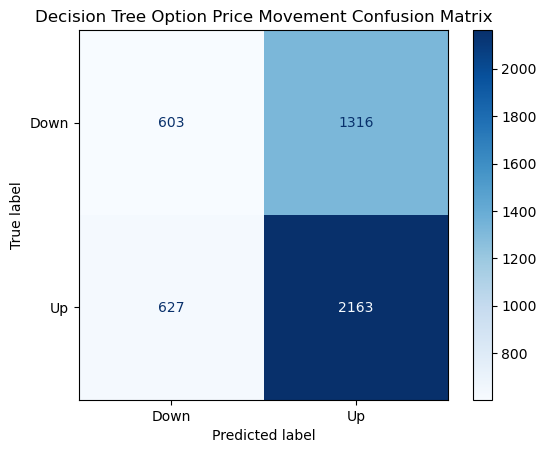}
    \caption{Confusion Matrix for option price predictions}
    \label{fig:tree2_matrix}
\end{figure}

\subsection{Decision Tree with Gradient Boosting}
\begin{figure}[ht!]
    \centering
    \includegraphics[width=1\textwidth]{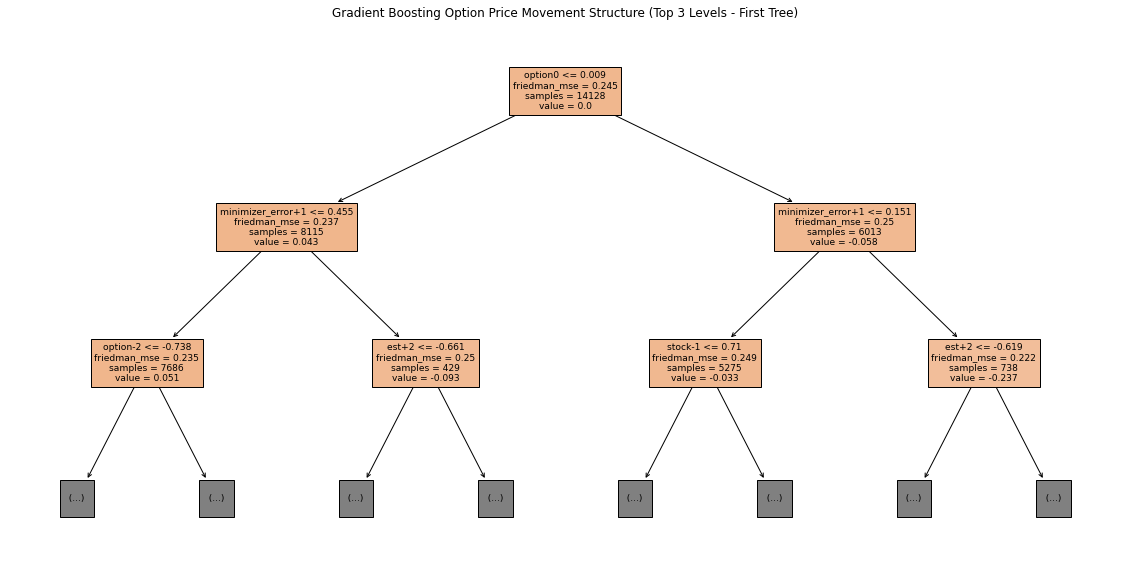}
    \caption{Gradient Boost Tree for option movement}
    \label{fig:tree4}
\end{figure}

Since the tree structure is too complex for demonstration purposes, we will only show a partial tree representing the first few layers \ref{fig:tree4}.
\begin{figure}[ht!]
    \centering
    \includegraphics[width=0.7\textwidth]{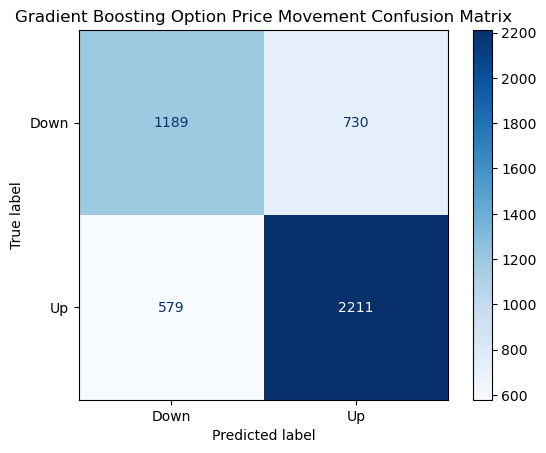}
    \caption{Confusion Matrix for option price predictions}
    \label{fig:tree4_matrix}
\end{figure}

\subsection{Random Forest}
\begin{figure}[ht!]
    \centering
    \includegraphics[width=1\textwidth]{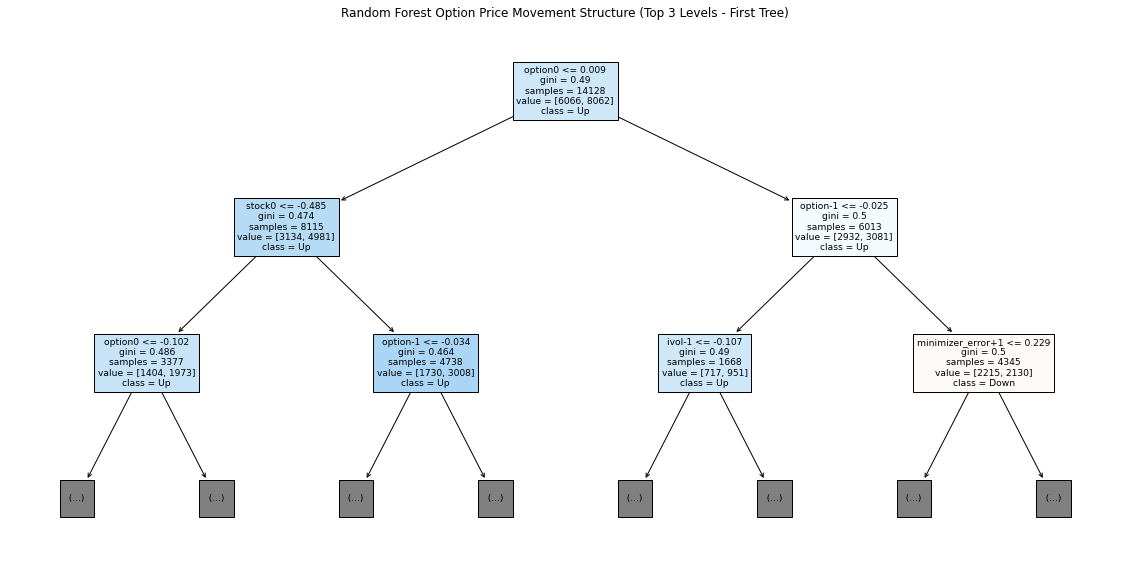}
    \caption{Random Forest Tree for option movement}
    \label{fig:tree6}
\end{figure}

Since the tree structure is too complex for demonstration purposes, we will only show a partial tree representing the first few layers \ref{fig:tree6}.
\begin{figure}[ht!]
    \centering
    \includegraphics[width=0.7\textwidth]{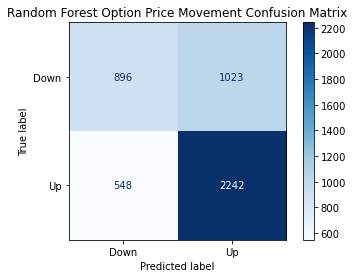}
    \caption{Confusion Matrix for option price predictions}
    \label{fig:tree6_matrix}
\end{figure}

\subsection{KNN}
\begin{figure}[ht!]
    \centering
    \includegraphics[width=0.75\textwidth]{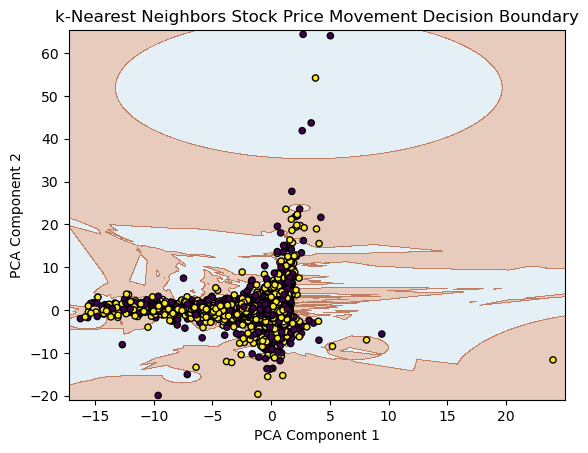}
    \caption{KNN boundary plot for option movement}
    \label{fig:knn}
\end{figure}

Since the KNN boundary is too complex for demonstration purposes, we will only show a partial boundary representing the first few boundaries \ref{fig:knn}.
\begin{figure}[ht!]
    \centering
    \includegraphics[width=0.7\textwidth]{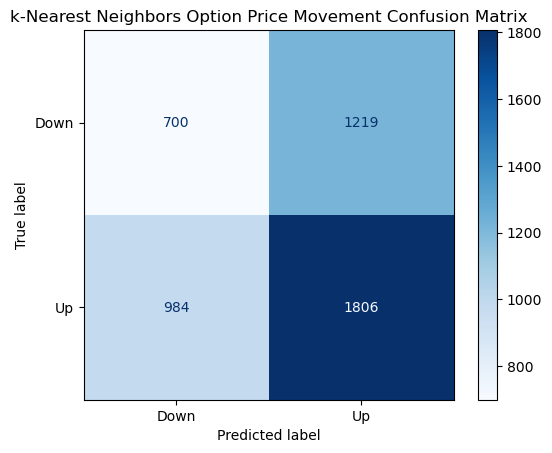}
    \caption{Confusion Matrix for option price predictions}
    \label{fig:knn2_matrix}
\end{figure}

\subsection{Neural Network}

\begin{figure}[!ht]
    \centering
    \includegraphics[width=0.77\textwidth]{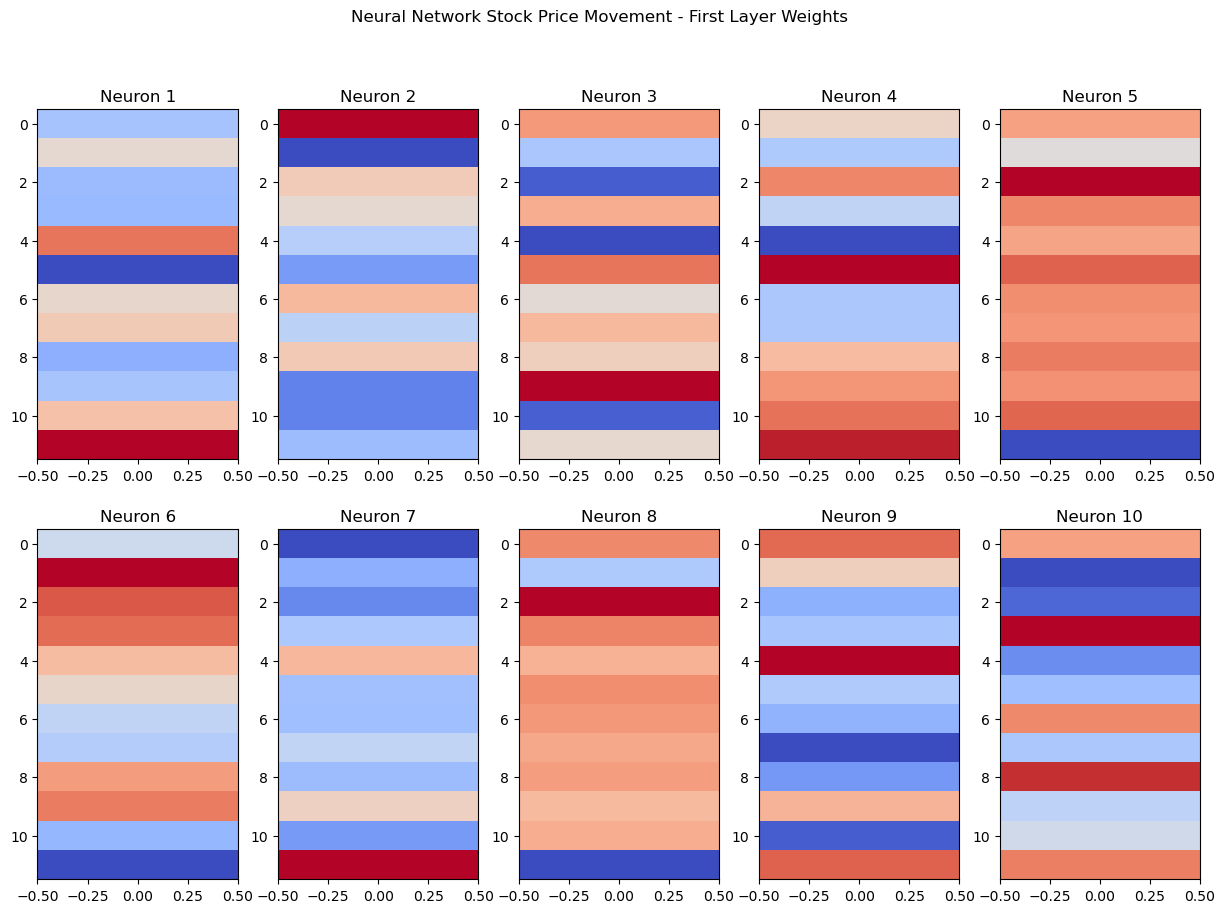}
    \caption{Neuron plot for option movement}
    \label{fig:nn}
\end{figure}

Since the Neural Network full layer is too complex for demonstration purposes, we will only show partial neuron graphs \ref{fig:nn}.
\begin{figure}[!ht]
    \centering
    \includegraphics[width=0.7\textwidth]{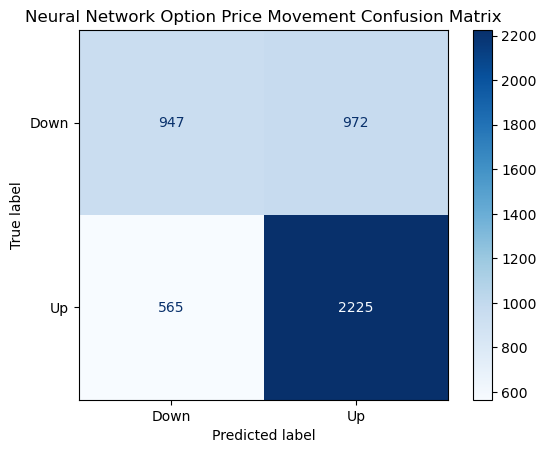}
    \caption{Confusion Matrix for option price predictions}
    \label{fig:nn2_matrix}
\end{figure}

\subsection{Stock Market Data:}
\begin{figure} 
    \centering
    \includegraphics[width=1\textwidth]{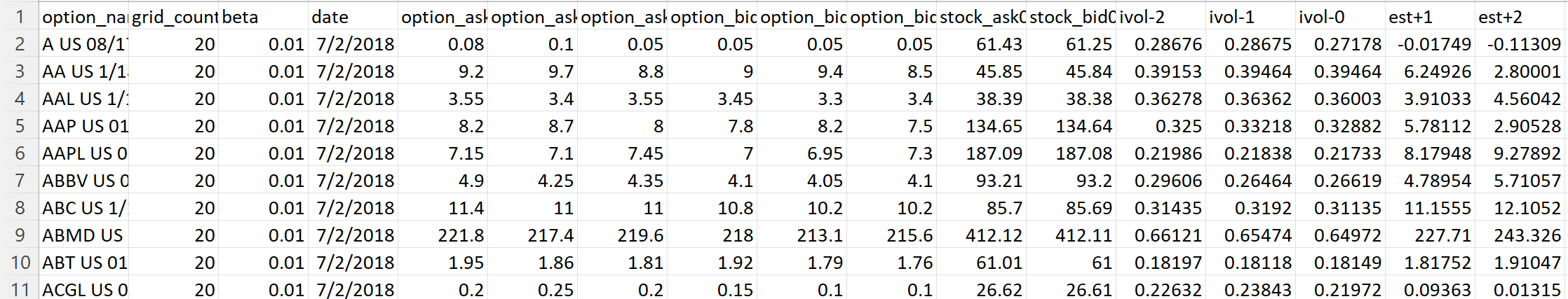}
    \caption{Real stock market data from Bloomberg.com}
    \label{fig:stock_data}
\end{figure}
The data used for this project were sourced from Bloomberg.com \cite{Bloom}. The dataset comprises stock prices, implied volatility, option prices, QRM estimators, and estimator errors across five different days for 23,548 unique stocks. This extensive dataset was selected to ensure a comprehensive representation of the entire stock market, avoiding biases toward any specific sector. We believe this diverse dataset provides a robust foundation for achieving more reliable and promising results. Additionally, the chosen features are instrumental in predicting future movements in options, which are subsequently utilized as inputs for our classification models.

\subsection{Result}
\begin{center}
\begin{table}
\caption{Prediction results by all classification models}
\label{table:results}
\begin{center}
\begin{tabular}{c| c| c| c} 
 \hline
 Method & Accuracy & Precision & Recall\\ [0.5ex] 
 \hline
  Decision Tree & 58.74\% & 62.17\% & 77.53\% \\
 \hline
  Gradient Boosting & 72.2\% & 75.18\% & 79.25\% \\
 \hline
 Random Forest & 66.64\% & 68.67\% & 80.38\% \\ [1ex] 
  \hline
 K Nearest Neighbors & 53.22\% & 59.7\% & 64.73\% \\ [1ex] 
  \hline
 Neural Network & 67.36\% & 69.6\% & 79.75\% \\ [1ex] 
 \hline
\end{tabular}
\end{center}
\end{table}
\end{center}

\section{Summary and Future Development}\label{sec:conclusions}
Although accurately predicting exact option prices remains challenging due to market volatility, we can make reliable predictions about whether an option's price will increase or decrease. These predictions are derived from market features extracted from the dataset of 23,548 entries and the minimizer obtained using the Quasi-Reversibility Method (QRM). Various classification models were tested, and their results are summarized in the table above. Based on our assumptions, we interpret the precision percentage as a key indicator of profitability when buying a call option. Using the Decision Tree model, we achieved a precision of 62.17$\%$, which improved to 68.67$\%$ with the Random Forest model and further to 75.18$\%$ when employing Gradient Boosting \ref{table:results}. Compared to prior results from Table 1 \ref{table:prev_results}, which showed a precision of 60.32$\%$ using a Neural Network, our current models demonstrate a notable improvement in precision and profitability potential.

However, the current approach only predicts whether a option will increase or decrease in value one day in advance, without accounting for factors such as transaction fees, holding periods, or longer-term market volatility. Additionally, the dataset spans a relatively short timeframe, which may limit its ability to capture longer-term market trends. Furthermore, the model currently classifies non-profitable options as price increases, which could skew the interpretation of results.

In future work, we plan to enhance our methodology by predicting the percentage change in option prices rather than simply classifying price movements as increases or decreases. This improvement will provide a more detailed understanding of market dynamics and help traders identify options with higher profit potential. Additionally, we aim to refine the Quasi-Reversibility Method (QRM) minimizer to improve its robustness and accuracy in handling noisy market data, incorporating techniques such as adaptive regularization for more reliable predictions. We also plan to explore advanced machine learning models beyond tree-based classifiers that are capable of capturing complex temporal patterns. 

\newpage
Appendix
\begin{center} 
\end{center}
First, we have some input for our classification:
\begin{align}
    \text{Data}(x_i,y_i)_{i=1}^{n} \tag{27}
\end{align}
Also, we have a loss function:
\begin{align}
    L(y_i, F(x)) \tag{28}
\end{align}
where
\begin{enumerate}
    \item $x_i$ - The input data fed into the model,
    \item $y_i$ - The output we are trying to predict.
\end{enumerate}
\textbf{Likelihood Function}: The likelihood function describes the joint probability of the observed data as a function of the parameters of our input:
\begin{align}
    L(\theta) = L(\theta; y) = f_y(y; \theta) \tag{29}
\end{align}
The likelihood function is the product of all such probabilities:
\begin{align}
L(\theta) = \prod_{i=1}^{n} f_i(y_i; \theta) \tag{30}
\end{align}
Here, we use the log-likelihood function:
\begin{align}
\log\_likelihood = y_i \cdot \log(p) + (1 - y_i) \cdot \log(1-p) \tag{31}
\end{align}
Here, $y_i$ must be binary (either 0 or 1), and $p$ represents the predicted probability.
Since we aim to improve predictions, we want to maximize the log-likelihood function. To use this as our loss function, we simply invert it:
\begin{align}
    \log\_loss = -\log\_likelihood \tag{32}
\end{align}
Minimizing the log-loss function corresponds to a better fit. Since the likelihood function is a probability function, we convert it into a function of corresponding odds using the logit function:
\begin{align}
\text{logit}(p) = \log\left(\frac{p}{1-p}\right) \tag{33}
\end{align}
where $\frac{p}{1-p}$ represents the odds of the given probability. Let's perform some algebraic manipulations:
\begin{align*}
    \log\_loss  &= - \left[ y_i \cdot \log\left(\frac{p}{1-p}\right) + \log(1-p) \right] \\
              &= - \left[ y_i \cdot \log(p) - y_i \cdot \log(1-p) + \log(1-p) \right] \\
              &= -y_i \cdot \log\left(\frac{p}{1-p}\right) - \log(1-p) \tag{34}
\end{align*}
Let $\log\left(\frac{p}{1-p}\right)$ be denoted as $y$. Then, solve for $p$:
\begin{align*}
    y &= \log\left(\frac{p}{1-p}\right) \\
    e^y &= \frac{p}{1-p} \\
    e^y + 1 &= \frac{1}{1-p} \\
    1-p &= \frac{1}{1+e^y} \\
    p &= 1 - \frac{1}{1+e^y} \\
    p &= \frac{e^y}{1+e^y} \tag{35}
\end{align*}
Substituting back:
\begin{align*}
    \log(1-p) &= \log\left(\frac{1}{1+e^y}\right) \\
              &= -\log(1+e^y) \tag{36}
\end{align*}
Substituting these into the log-loss function:
\begin{align}
    \log\_loss &= -y_i \cdot \log\left(\frac{p}{1-p}\right) + \log(1+e^y) \tag{37}
\end{align}
Differentiating the loss function:
\begin{align}
    \frac{d}{d \log\left(\frac{p}{1-p}\right)}\left(-y_i \cdot \log\left(\frac{p}{1-p}\right) + \log(1+e^{\log(\frac{p}{1-p})})\right) = -y_i + \frac{e^{\log(\frac{p}{1-p})}}{1+e^{\log(\frac{p}{1-p})}} \tag{38}
\end{align}
This can be simplified to, Predicted minus Observed:
\begin{align}
    \text{Observed} = y_i \tag{39} \\
    \text{Predicted} = \frac{e^{\log(\frac{p}{1-p})}}{1+e^{\log(\frac{p}{1-p})}} \tag{40}
\end{align}
\cite{books}

\newpage
%
%

\end{document}